 \documentclass[11pt]{article}

 \usepackage{mathrsfs,amsfonts,amsmath,amssymb}
 \usepackage[dvips]{color}

 \newtheorem{definition}{Definition}[section]
 \newtheorem{hypothesis}{Hypothesis}[section]
 \newtheorem{lemma}{Lemma}[section]
 \newtheorem{proposition}{Proposition}[section]
 \newtheorem{theorem}{Theorem}[section]
 \newtheorem{corollary}{Corollary}[section]
 \newtheorem{example}{Example}[section]

 \def\blemma{\begin{lemma}\sl{}\def\elemma{\end{lemma}}}
 \def\bproposition{\begin{proposition}\sl{}\def\eproposition{\end{proposition}}}
 \def\btheorem{\begin{theorem}\sl{}\def\etheorem{\end{theorem}}}
 
 \def\bexample{\begin{example}\rm{}\def\eexample{\end{example}}}

 \def\beqlb{\begin{eqnarray}}\def\eeqlb{\end{eqnarray}}
 \def\beqnn{\begin{eqnarray*}}\def\eeqnn{\end{eqnarray*}}

 \def\mbb{\mathbb}
 \def\mbf{\mathbf}
 \def\mcr{\mathscr}

 \def\proof{\noindent\textit{Proof.~}}\def\qed{\hfill$\Box$\medskip}

 \def\<{\langle}\def\>{\rangle}

 \def\nnm{\nonumber}

\begin{document}

\noindent{\small Published in: \textit{Frontiers of Mathematics in
China} \textbf{1} (2006), 1: 73--97.}

\bigskip\bigskip

\centerline{\LARGE\textbf{Branching processes with immigration}}

\smallskip

\centerline{\LARGE\textbf{and related topics}}

\bigskip

\centerline{Zenghu LI\,\footnote{~E-mail: lizh@bnu.edu.cn}}

\medskip

\centerline{\small School of Mathematical Sciences, Beijing Normal
University, Beijing 100875, China}

\bigskip\bigskip

{\narrower{\narrower

\noindent\textbf{Abstract:} This is a survey on recent progresses
in the study of branching processes with immigration, generalized
Ornstein-Uhlenbeck processes and affine Markov processes. We
mainly focus on the applications of skew convolution semigroups
and the connections in those processes.

\medskip


\medskip


\par}\par}

\bigskip

\section{Introduction}

\setcounter{equation}{0}

Let $\mbb{N} := \{0,1,2,\cdots\}$ and let $\{\xi(n,i): n,i = 1,2,
\cdots\}$ be a sequence of $\mbb{N}$-valued i.i.d.\ random
variables. Let $x(0)$ be an $\mbb{N}$-valued random variable which
is independent of $\{\xi(n,i)\}$. A \textit{Galton-Watson
branching process} (GW-process) $\{x(n): n = 0,1,2, \cdots\}$ is
defined inductively by
 \beqlb\label{1.1}
x(n) = \sum_{i=1}^{x(n-1)} \xi(n,i), \qquad n=1,2, \cdots.
 \eeqlb
This process is a mathematical representation of the random
evolution of an isolated population. We refer the reader to
Athreya and Ney \cite{AN72} and Harris \cite{Har63} for the theory
of branching processes.

A useful and realistic modification of the above scheme is the
addition of the possibility of immigration into the population.
From the point of applications, the immigration processes are
clearly of great importance. Let $\{\eta(n): n = 1,2, \cdots\}$ be
another sequence of $\mbb{N}$-valued i.i.d.\ random variables
which are independent of $\{\xi(n,i)\}$. Let $y(0)$ be an
$\mbb{N}$-valued random variable independent of $\{\xi(n,i)\}$ and
$\{\eta(n)\}$. We can define a \textit{Galton-Watson branching
process with immigration} (GWI-process) $\{y(n): n =
0,1,2,\cdots\}$ by
 \beqlb\label{1.2}
y(n) = \sum_{i=1}^{y(n-1)} \xi(n,i) + \eta(n), \qquad n=1,2,
\cdots;
 \eeqlb
see, e.g., \cite[p.263]{AN72}. The intuitive meaning of the
process is clear from the construction \eqref{1.2}. Let $g(\cdot)$
and $h(\cdot)$ be the generating function of $\{\xi(n,i)\}$ and
$\{\eta(n)\}$, respectively. Because of the independence of the
random variables $\{\xi(n,i), \eta(n): n,i = 1,2, \cdots\}$ it is
easy to see that $\{y(n)\}$ is a discrete-time Markov chain with
one-step transition matrix $P(i,j)$ defined by
 \beqlb\label{1.3}
\sum_{j=0}^\infty P(i,j)z^j = g(z)^ih(z), \qquad 0\le z\le 1; i=
0,1,2,\cdots.
 \eeqlb

The purpose of this survey is to give a brief introduction to the
recent progresses in the study of branching processes with
immigration and related topics. We shall be concerned with
continuous state branching processes (CB processes), CB processes
with immigration (CBI processes), measure-valued branching
processes (MB processes), Dawson-Watanabe superprocesses,
immigration superprocesses, generalized Ornstein-Uhlenbeck
processes, and affine processes. The basic mathematical structures
of those processes already exist in \eqref{1.1}--\eqref{1.3}. Our
emphasis is on the applications of skew convolution semigroups and
the connections in those processes. This is an ongoing research
topic of the Probability Group in Beijing Normal University. Other
topics where our Group has been involved include interacting
particle systems, ergodic and spectral theory, probabilistic and
functional inequalities, large and moderate deviations and so on.
We refer the reader to Chen \cite{Che02, Che04a, Che04b} and Wang
\cite{Wan05a, Wan05b, Wan05c} for some recent results of the Group
on those topics.

Let us introduce some notation which will be used throughout the
survey. Given a metrizable topological space $E$, we denote by
$\mcr{B}(E)$ its Borel $\sigma$-algebra. Let $B(E)$ the space of
bounded real $\mcr{B}(E)$-measurable functions on $E$ and $C(E)$
the subset of $B(E)$ of continuous functions. Let $M(E)$ be the
space of finite Borel measures on $E$ endowed with the topology of
weak convergence. For $f\in B(E)$ and $\mu\in M(E)$, let $\mu(f) =
\int_E fd\mu$. Let $\delta_x$ denote the unit mass concentrated at
$x\in E$. For any integer $m\ge 1$ let $C^m(\mbb{R}^d)$ denote the
set of smooth functions on the Euclidean space $\mbb{R}^d$ with
all partial derivatives up to the $m$th order belonging to
$C(\mbb{R}^d)$. Let $C^\infty(\mbb{R}^d) = \bigcap_{m=1}^\infty
C^m(\mbb{R}^d)$.

\section{Skew convolution semigroups and examples}

\setcounter{equation}{0}

Let $(S,+)$ be a metrizable abelian semigroup, that is, $S$ is a
metrizable topological space and there is a composition law $+:
S^2 \to S$ which is associative, commutative and continuous. For
two Borel probability measures $\mu$ and $\nu$ on $S$, the image
of the product measure $\mu \times \nu$ under the composition law
is called the \textit{convolution} of $\mu$ and $\nu$ and is
denoted by $\mu*\nu$. Suppose that $(Q_t)_{t\ge0}$ is a Borel
Markov transition semigroup on $S$ satisfying $Q_t(0,\cdot) =
\delta_0$ and the \textit{branching property}
 \beqlb\label{2.1}
Q_t(x_1+x_2,\cdot) = Q_t(x_1,\cdot) * Q_t(x_2,\cdot), \qquad
t\ge0, x_1,x_2\in S.
 \eeqlb
Given $t\ge0$ and a Borel measure $\mu$ on $S$ we define the
measure $\mu Q_t$ by
 \beqnn
\mu Q_t(A) = \int_S Q_t(x,A) \mu(dx), \qquad A\in \mcr{B}(S).
 \eeqnn

\blemma\label{l2.1} For any Borel probability measures $\mu$ and
$\nu$ on $S$ we have
 \beqlb\label{2.2}
(\mu*\nu)Q_t = (\mu Q_t)*(\nu Q_t), \qquad t\ge0.
 \eeqlb
\elemma

\proof Let $f\in B(S)$. From the branching property it follows
that
 \beqnn
\int_S Q_tf(x)(\mu*\nu)(dx)
 &=&
\int_S \mu(dx)\int_S Q_tf(x+y)\nu(dy)  \\
 &=&
\int_S \mu(dx)\int_S \nu(dy)\int_S f(z)Q_t(x+y,dz)  \\
 &=&
\int_S \mu(dx)\int_S \nu(dy)\int_SQ_t(x,dz_1)\int_S f(z_1+z_2)
Q_t(y,dz_2)  \\
 &=&
\int_S (\mu Q_t)(dz_1)\int_S f(z_1+z_2)(\nu Q_t)(dz_2)  \\
 &=&
\int_S f(z)[(\mu Q_t)*(\nu Q_t)](dz).
 \eeqnn
Then we have the equality \eqref{2.2}. \qed

\btheorem\label{t2.1} Suppose that $(\gamma_t)_{t\ge0}$ is a
family of Borel probability measures on $S$. Then
 \beqlb\label{2.3}
Q^\gamma_t (x,\cdot) := Q_t(x,\cdot) * \gamma_t(\cdot), \qquad
x\in S, t\ge0
 \eeqlb
defines a Borel kernel on $S$ and $(Q_t^\gamma)_{t\ge0}$ form a
transition semigroup if and only if
 \beqlb\label{2.4}
\gamma_{r+t} = (\gamma_rQ_t) * \gamma_t, \qquad r,t\ge0.
 \eeqlb
\etheorem

\proof It is easy to show that $Q_t^\gamma(x,dy)$ is a Borel
kernel on $S$. Then we only need to prove that \eqref{2.4} is
equivalent to the Chapman-Kolmogorov equation
 \beqlb\label{2.5}
\int_S f(y)Q_{r+t}^\gamma(x,dy) = \int_SQ_r^\gamma(x,dy)\int_S
f(z)Q_t^\gamma(y,dz), \qquad f\in B(S).
 \eeqlb
If \eqref{2.5} holds, we may apply this equation with $x=0$ to see
that
 \beqnn
\int_S f(z)\gamma_{r+t}(dz)
 &=&
\int_S \gamma_r(dy)\int_S f(z)Q_t^\gamma(y,dz)  \\
 &=&
\int_S \gamma_r(dy)\int_SQ_t(y,dz_1)\int_S f(z_1+z_2)\gamma_t(dz_2)  \\
 &=&
\int_S (\gamma_rQ_t)(dz_1)\int_S f(z_1+z_2)\gamma_t(dz_2).
 \eeqnn
Then \eqref{2.4} holds. Conversely, if \eqref{2.4} holds, we have
 \beqnn
\int_S f(z)Q_{r+t}^\gamma(x,dz)
 &=&
\int_S Q_{r+t}(x,dz_1)\int_S f(z_1+z_2)\gamma_{r+t}(dz_2)  \\
 &=&
\int_S Q_r(x,dy)\int_SQ_t(y,dz_1)\int_S (\gamma_rQ_t)(dz_2)
\int_S f(z_1+z_2+z_3)\gamma_t(dz_3)  \\
 &=&
\int_{S}Q^\gamma_r(x,dy)\int_SQ_t(y,dz_2)
\int_S f(z_2+z_3)\gamma_t(dz_3)  \\
 &=&
\int_SQ_r^\gamma(x,dy)\int_S f(z)Q_t^\gamma(y,dz).
 \eeqnn
That proves the Chapman-Kolmogorov equation \eqref{2.5}. \qed

We call $(\gamma_t)_{t\ge0}$ a \textit{skew convolution semigroup}
(SC-semigroup) associated with $(Q_t)_{t\ge0}$ if it satisfies
\eqref{2.4}; see Li \cite{Li95/6, Li02}. The the kernels
$Q_t^\gamma(x,dy)$ defined by \eqref{2.3} give an abstract
formulation of the expression \eqref{1.3}. In particular, if $Q_t$
is the identity operator for every $t\ge0$, the SC-semigroup
defined by \eqref{2.4} becomes a standard convolution semigroup
and $(Q^\gamma_t)_{t\ge0}$ is the transition semigroup of a
\textit{L\'evy process}. We refer the reader to Bertoin
\cite{Ber96} and Sato \cite{Sat99} for the theory of L\'evy
processes. The general formulae \eqref{2.3} and \eqref{2.4}
include many additional mathematical contents, which are
illustrated by the following examples.

\bexample\label{e1.1} In the particular case $S = \mbb{R}_+$, a
Markov process with transition semigroup $(Q_t)_{t\ge0}$ is called
a \textit{CB-process} and a Markov process with transition
semigroup $(Q_t^\gamma)_{t\ge0}$ is called a \textit{CBI-process};
see \cite{KW71, SW73}. \eexample

\bexample\label{e1.2} If $S = M(E)$ is the space of all finite
Borel measures on a metrizable space $E$, the semigroup
$(Q_t)_{t\ge0}$ corresponds to an MB-process, of which the
\textit{Dawson-Watanabe superprocess} is a special case; see
\cite{Daw93}. A Markov process with state space $M(E)$ is
naturally called an \textit{immigration superprocess} associated
with $(Q_t)_{t\ge0}$ if it has transition semigroup
$(Q_t^\gamma)_{t\ge0}$; see \cite{Li95/6, Li96, Li02}. \eexample

\bexample\label{e1.3} Let us consider the case where $S=H$ is a
real separable Hilbert space and $Q_t(x,\cdot) \equiv
\delta_{T_tx}$ for a strongly continuous semigroup of bounded
linear operators $(T_t)_{t\ge0}$ on $H$. In this case,
$(Q^\gamma_t)_{t\ge0}$ is called a \textit{generalized Mehler
semigroup} associated with $(T_t)_{t\ge0}$, which corresponds to a
generalized Ornstein-Uhlenbeck process (OU-process). This
formulation of the processes was given by Bogachev \textit{et
al.}\ \cite{BRS96}; see also \cite{DLSS04, FR00}. \eexample

\bexample\label{e1.4} If $S = \mbb{R}_+^m \times \mbb{R}^n$ for
integers $m\ge0$ and $n\ge0$, the transition semigroup
$(Q^\gamma_t)_{t\ge0}$ corresponds to an \textit{affine process}.
The affine Markov processes were introduced in mathematical
finance; see, e.g., \cite{DL06, DFS03}. \eexample

\section{Continuous state branching processes with immigration}

\setcounter{equation}{0}

There is a rich literature in the study of CB- and CBI-processes.
In particular, the class of CBI-processes was characterized
completely by Kawazu and Watanabe \cite{KW71}. Let $F$ be a
function defined by
 \beqlb\label{3.1}
F(\lambda) = b\lambda + \int_0^\infty(1-e^{-\lambda u})m(du),
\qquad \lambda\ge 0,
 \eeqlb
where $b\ge0$ is a constant and $um(du)$ is a finite measure on
$(0,\infty)$. Let $R$ be given by
 \beqlb\label{3.2}
R(\lambda) = \beta\lambda - \alpha\lambda^2 - \int_0^\infty
\big(e^{-\lambda u} - 1 + \lambda u\big)\mu(du), \qquad \lambda\ge
0,
 \eeqlb
where $\beta \in \mbb{R}$ and $\alpha\ge0$ are constants and
$(u\land u^2)\mu(du)$ is a finite measure on $(0,\infty)$. We can
define a transition semigroup $(P_t)_{t\ge0}$ on $\mbb{R}_+$ by
 \beqlb\label{3.3}
\int_0^\infty e^{-\lambda y} P_t(x,dy)
 =
\exp\bigg\{-x\psi_t(\lambda) - \int_0^t F(\psi_s(\lambda))
ds\bigg\}, \qquad \lambda\ge 0,
 \eeqlb
where $\psi_t(\lambda)$ is the unique solution of
 \beqlb\label{3.4}
\frac{d\psi_t}{dt}(\lambda) = R(\psi_t(\lambda)), \qquad
\psi_0(\lambda) = \lambda.
 \eeqlb
A Markov process $\{y(t): t\ge0\}$ with transition semigroup
$(P_t)_{t\ge0}$ is a special case of the CBI-process defined in
\cite{KW71}. Here \eqref{3.3} is the continuous time version of
\eqref{1.3}. Various limit theorems for the CBI-process have been
established; see, e.g., \cite{Gre74, Li00a, Pak88, Pak99} and the
references therein.

The connections between the GWI-processes and the CBI-processes
were investigated in Kawazu and Watanabe \cite{KW71}. They showed
that a CBI-processes arises as the high density limit in
finite-dimensional distributions of a sequence of GWI-processes.
Some simple conditions were given in Li \cite{Li05} which ensure
that the convergence of GWI-processes mentioned above holds on the
space of c\`adl\`ag paths. Let $\{y_k(n): n\ge0\}$ be a sequence
of GWI-processes with parameters $\{(g_k,h_k)\}$ and
$\{\gamma_k\}$ a sequence of positive numbers. For $0\le
\lambda\le k$ set
 \beqlb\label{3.5}
F_k(\lambda) = \gamma_k[1 - h_k(1-\lambda/k)]
 \eeqlb
and
 \beqlb\label{3.6}
R_k(\lambda) = k\gamma_k[(1-\lambda/k) - g_k(1-\lambda/k)].
 \eeqlb
Let us consider the following conditions: \begin{itemize}

\item[(3.A)]
As $k\to \infty$, we have $\gamma_k \to\infty$ and $\gamma_k/k
\to$ some $\gamma_0\ge 0$.

\item[(3.B)]
The sequence $\{F_k\}$ defined by \eqref{3.5} is uniformly
Lipschitz on each bounded interval and converges as $k\to \infty$.

\item[(3.C)]
The sequence $\{R_k\}$ defined by \eqref{3.6} is uniformly
Lipschitz on each bounded interval and converges as $k\to \infty$.
\end{itemize}

\btheorem\label{t3.1} {\rm (\cite{Li05})} Suppose that conditions
(3.A), (3.B) and (3.C) are satisfied. If $y_k(0)/k$ converges in
distribution to $y(0)$, then $\{y_k([\gamma_kt])/k: t\ge0\}$
converges in distribution on $D([0,\infty), \mbb{R}_+)$ to the
CBI-process $\{y(t): t\ge0\}$ corresponding to $(R,F)$. \etheorem

\proof We here only give a sketch and refer the reader to
\cite{Li05} for the details. Let $A$ denote the generator of the
CBI-process. For $\lambda>0$ and $x\ge0$ set $e_\lambda(x) =
e^{-\lambda x}$ and let $D$ be the linear hull of $\{e_\lambda:
\lambda>0\}$. For $\lambda>0$ we have
 \beqlb\label{3.7}
Ae_\lambda(x)
 =
-e^{-\lambda x}\left[xR(\lambda) + F(\lambda)\right], \qquad x\in
\mbb{R}_+,
 \eeqlb
and this equality determines the actions of $A$ on $D$ by
linearity. Then we deduce that $D$ is a core of $A$. Note that
$\{y_k (n)/k: n\ge0\}$ is a Markov chain with state space $E_k :=
\{0,1/k,2/k,\cdots\}$ and one-step transition probability
$Q_k(x,dy)$ determined by
 \beqnn
\int_{E_k}e^{-\lambda y}Q_k(x,dy)
 =
g_k(e^{-\lambda/k})^{kx}h_k(e^{-\lambda/k}).
 \eeqnn
Then one checks that the (discrete) generator $A_k$ of $\{y_k
([\gamma_kt])/k: t\ge0\}$ is given by
 \beqnn
A_ke_\lambda(x)
 &=&
\gamma_k\Big[g_k(e^{-\lambda/k})^{kx}h_k(e^{-\lambda/k})
- e^{-\lambda x}\Big] \\
 &=&
\gamma_k\Big[\exp\{xk\alpha_k(\lambda)(g_k(e^{-\lambda/k})-1)\}
\exp\{\beta_k(\lambda)(h_k(e^{-\lambda/k})-1)\} - e^{-\lambda
x}\Big],
 \eeqnn
where
 \beqnn
\alpha_k(\lambda)
 =
(g_k(e^{-\lambda/k})-1)^{-1} \log g_k(e^{-\lambda/k})
 \eeqnn
and $\beta_k(\lambda)$ is defined by the same formula with $g_k$
replaced by $h_k$. Then we use the assumptions to show that
 \beqlb\label{3.8}
A_ke_\lambda(x)
 =
-e^{-\lambda x}\left[x\alpha_k(\lambda)S_k(\lambda) +
x\gamma_k(\alpha_k(\lambda)-1)\lambda + H_k(\lambda)\right] +
o(1),
 \eeqlb
where
 \beqnn
H_k(\lambda)
 =
\gamma_k\beta_k(\lambda)(1-h_k(e^{-\lambda/k})).
 \eeqnn
By elementary calculations we find that
 \beqnn
\alpha_k(\lambda)
 =
1 + \frac{1}{2}(1-g_k(e^{-\lambda/k})) + o(1-g_k(e^{-\lambda/k})),
 \eeqnn
and so $\lim_{k\to\infty} \gamma_k(\alpha_k(\lambda)-1) =
\gamma_0\lambda/2$. It follows that
 \beqnn
\lim_{k\to\infty}\left[\alpha_k(\lambda)S_k(\lambda) +
\gamma_k(\alpha_k(\lambda)-1)\lambda\right] = R(\lambda).
 \eeqnn
Then one shows that $\lim_{k\to\infty}H_k(\lambda) =
\lim_{k\to\infty} F_k(\lambda) = F(\lambda)$. In view of
\eqref{3.7} and \eqref{3.8} we get
 \beqnn
\lim_{k\to\infty} \sup_{x\in E_k}\left|A_ke_\lambda(x) -
Ae_\lambda(x)\right| = 0
 \eeqnn
for each $\lambda>0$. This clearly implies that
 \beqnn
\lim_{k\to\infty} \sup_{x\in E_k}\left|A_kf(x) - Af(x)\right| = 0
 \eeqnn
for each $f\in D$. That proves the desired convergence. \qed

By the results of Li \cite{Li91} it is easy to show that the limit
functions of $\{F_k\}$ and $\{R_k\}$ always have representations
\eqref{3.5} and \eqref{3.6}, respectively. On the other hand, for
any $(F,R)$ given by \eqref{3.1} and \eqref{3.2}, there are
sequences $\{\gamma_k\}$ and $\{(g_k, h_k)\}$ as above such that
(3.A), (3.B) and (3.C) hold with $F_k \to F$ and $R_k \to R$; see
\cite{Li92, Li05}. Those results show the range of applications of
Theorem~\ref{t3.1}. As consequences of the above theorem, Li
\cite{Li05} gave some generalizations of the Ray-Knight Theorems
on Brownian local times; see also Le Gall and Le Jan \cite{LL98}.
We remark that conditions (3.A), (3.B) and (3.C) parallel the
sufficient conditions for the convergence of continuous-time and
discrete state branching processes with immigration, see, e.g.,
\cite{Li92}. In most cases, those conditions are easier to check
than the sufficient conditions given by Kawazu and Watanabe
\cite{KW71}, which involve complicated composition and convolution
operations.

From \eqref{3.3} it is easy to see that the transition semigroup
$(P_t)_{t\ge0}$ is Fellerian, so the CBI-process has a Hunt
realization. A construction of the process was given in Dawson and
Li \cite{DL06} as the strong solution of a stochastic integral
equation. Suppose that $(\Omega, \mcr{F}, \mcr{F}_t, \mbf{P})$ is
a filtered probability space satisfying the usual hypotheses on
which the following adapted objects are defined:
\begin{itemize}

\item
a standard Brownian motion $\{B(t)\}$;

\item
a Poisson random measure $N_0(ds,d\xi)$ on $\mbb{R}_+^2$ with
intensity $dsm(d\xi)$;

\item
a Poisson random measure $N_1(ds,du,d\xi)$ on $\mbb{R}_+^3$ with
intensity $dsdu\mu(d\xi)$;

\end{itemize}
We assume that $\{B(t)\}$, $\{N_0(ds,d\xi)\}$ and
$\{N_1(ds,du,d\xi)\}$ are independent of each other. Let $x(0)$ be
a non-negative $\mcr{F}_0$-measurable random variable satisfying
$\mbf{E}[x(0)] < \infty$. We consider the stochastic integral
equation
 \beqlb\label{3.9}
x(t) &=& x(0) + \int_0^t(b + \beta x(s))ds +
\int_0^t \sqrt{2\alpha x(s)} dB(s) \nnm \\
 & &
+ \int_0^t\int_0^\infty \xi N_0(ds,d\xi) + \int_0^t \int_0
^{x(s-)}\int_0^\infty \xi \tilde N_1(ds,du,d\xi),
 \eeqlb
where $\tilde N_1(ds,du,d\xi) = N_1(ds,du,d\xi) - dsdu\mu(d\xi)$.

\btheorem\label{t3.2} {\rm (\cite{DL06})} There is a unique
non-negative c\`adl\`ag process $\{x(t): t\ge0\}$ such that
equation \eqref{3.9} is satisfied a.s.\ for every $t\ge0$.
\etheorem

The above theorem implies that \eqref{3.9} has a unique strong
solution $\{x(t): t\ge0\}$ and the solution is a strong Markov
process. For $f \in C^2(\mbb{R}_+)$ we see from \eqref{3.9} and
It\^o's formula (see, e.g., \cite[p.334-335]{DM82}) that
 \beqnn
f(x(t))
 &=&
f(x(0)) + \int_0^t f^\prime(x(s))(b+\beta x(s))ds +
\mbox{martingale} \nnm \\
 & &
+ \int_0^t\int_0^\infty f^\prime(x(s))\xi N_0(ds,d\xi)
+ \alpha\int_0^t f^{\prime\prime}(x(s)) x(s) ds \nnm \\
 & &
+ \int_0^t\int_0^\infty [f(x(s)+\xi) - f(x(s))
- f^\prime(x(s))\xi] N_0(ds,d\xi) \nnm \\
 & &
+ \int_0^t \int_0^{x(s-)}\int_0^\infty [f(x(s)+\xi) - f(x(s)) -
f^\prime(x(s))\xi] N_1(ds,du,d\xi) \nnm \\
 &=&
f(x(0)) + \int_0^t f^\prime(x(s))(b+\beta x(s))ds +
\mbox{martingale} \nnm \\
 & &
+\, \alpha\int_0^t f^{\prime\prime}(x(s)) x(s) ds
+ \int_0^tds\int_0^\infty [f(x(s)+\xi) - f(x(s))] m(d\xi) \nnm \\
 & &
+ \int_0^tds \int_0^\infty [f(x(s)+\xi) - f(x(s)) -
f^\prime(x(s))\xi] x(s) \mu(d\xi).
 \eeqnn
Then $\{x(t): t\ge0\}$ has generator $A$ defined by
 \beqlb\label{3.10}
A f(x)
 &=&
\alpha xf^{\prime\prime}(x) + (b + \beta x) f^\prime(x)
+ \int_0^\infty \big[f(x+\xi) - f(x)\big] m(d\xi) \nnm \\
 & &
+\, \int_0^\infty \big[f(x+\xi) - f(x) - f^\prime(x)\xi\big]x
\mu(d\xi),
 \eeqlb
so it is a CBI-process; see \cite{KW71}.

The approach of stochastic equations was also used in \cite{DL06}
to construct a general type of CBI-processes in random catalysts.
Suppose we have the parameters $(a, (\alpha_{ij}), (b_1,b_2),
(\beta_{ij}), m, \mu)$ such that
\begin{itemize}

\item
$a\in \mbb{R}_+$ is a constant;

\item
$(\alpha_{ij})$ is a symmetric non-negative definite $(2\times
2)$-matrix;

\item
$(b_1,b_2)\in \mbb{R}_+^2$ is a vector;

\item
$(\beta_{ij})$ is a $(2\times 2)$-matrix with $\beta_{12}=0$;

\item
$m(d\xi)$ is a $\sigma$-finite measure on $\mbb{R}_+^2$ supported
by $\mbb{R}_+^2 \setminus \{0\}$ such that
 \beqnn
\int_{\mbb{R}_+^2} [\xi_1 + \xi_2] m(d\xi) < \infty;
 \eeqnn

\item
$\mu(d\xi)$ is a $\sigma$-finite measure on $\mbb{R}_+^2$
supported by $\mbb{R}_+^2 \setminus \{0\}$ such that
 \beqnn
\int_{\mbb{R}_+^2} \big[(\xi_1\land \xi_1^2) + (\xi_2\land
\xi_2^2)\big] \mu(d\xi) < \infty.
 \eeqnn
\end{itemize}
Let $\sigma_0 = \sqrt{a}$ and let $(\sigma_{ij})$ be a $(2\times
2)$-matrix satisfying $(\alpha_{ij}) = (\sigma_{ij})
(\sigma_{ij})^\tau$. Let $(\Omega, \mcr{F}, \mcr{F}_t, \mbf{P})$
be a filtered probability space satisfying the usual hypotheses.
Suppose that on this probability space the following adapted
objects are defined:
\begin{itemize}

\item
a $3$-dimensional Brownian motion $\{(B_0(t), B_1(t), B_2(t))\}$;

\item
a Poisson random measure $N_0(ds,d\xi)$ on $\mbb{R}_+^3$ with
intensity $dsm(d\xi)$;

\item
a Poisson random measure $N_1(ds,du,d\xi)$ on $\mbb{R}_+^4$ with
intensity $dsdu\mu(d\xi)$.

\end{itemize}
We assume that $\{(B_0(t), B_1(t), B_2(t))\}$, $\{N_0(ds,d\xi)\}$
and $\{N_1(ds,du,d\xi)\}$ are independent of each other. Let
$x(0)$ and $y(0)$ be non-negative $\mcr{F}_0$-measurable random
variables defined on $(\Omega, \mcr{F}, \mcr{F}_t, \mbf{P})$. We
consider the equation system
 \beqlb
x(t) &=& x(0) + \int_0^t(b_1 + \beta_{11} x(s))ds +
\int_0^t\sigma_{11}\sqrt{2x(s)} dB_1(s) \nnm \\
 & &
+ \int_0^t\sigma_{12}\sqrt{2x(s)} dB_2(s) +
\int_0^t\int_{\mbb{R}_+^2} \xi_1 N_0(ds,d\xi) \nnm \\
 & &
+ \int_0^t \int_0 ^{x(s-)}\int_{\mbb{R}_+^2} \xi_1 \tilde
N_1(ds,du,d\xi), \label{3.11} \\
y(t) &=& y(0) + \int_0^t(b_2 + \beta_{21} x(s)y(s) + \beta_{22}
y(s))ds + \int_0^t\sigma_0\sqrt{2y(s)} dB_0(s) \nnm \\
 & &
+ \int_0^t\sigma_{21}\sqrt{2x(s)y(s)} dB_1(s)
+ \int_0^t\sigma_{22}\sqrt{2x(s)y(s)} dB_2(s)  \nnm \\
 & &
+ \int_0^t\int_{\mbb{R}_+^2} \xi_2 N_0(ds,d\xi) + \int_0^t
\int_0^{lx(s-)y(s-)} \int_{\mbb{R}_+^2} \xi_2 \tilde
N_1(ds,du,d\xi). \label{3.12}
 \eeqlb

\btheorem\label{t3.3} {\rm (\cite{DL06})} The equation system
given by \eqref{3.11} and \eqref{3.12} has a unique strong
solution $\{(x(t), y(t))\}$. \etheorem

Following Dawson and Fleischmann \cite{DF97}, we call $\{(x(t),
y(t))\}$ a \textit{catalytic CBI-process}, where $\{(x(t))\}$ is
the \textit{catalyst process} and $\{y(t)\}$ is the
\textit{reactant process}. It is not hard to see that $\{x(t)\}$
is a CBI-process.  Intuitively, we may think of $\{y(t)\}$ as a
CBI-process with random branching catalysts governed by the
process $\{x(t)\}$. A slightly more general catalytic CBI-process
with two reactant processes was considered in \cite{DL06}.

\section{Two-dimensional affine processes}

\setcounter{equation}{0}

The concept of affine Markov processes was introduced in the study
of financial models; see, e.g., \cite{DFS03} and the references
therein. For simplicity we only consider those processes in the
two-dimensional case. Let $D = \mbb{R}_+ \times \mbb{R}$ and $U =
\mbb{C}_- \times (i\mbb{R})$, where $\mbb{C}_- = \{a+ib: a\in
\mbb{R}_-, b\in \mbb{R}\}$ and $i\mbb{R} = \{ib: b\in \mbb{R}\}$.
A transition semigroup $(Q_t)_{t\ge0}$ on $D$ is called a
\textit{homogeneous affine semigroup} (HA-semigroup) if for each
$t\ge0$ there exists a continuous operator $u \mapsto \psi(t,u)$
on $U$ such that
 \beqlb\label{4.1}
\int_D \exp\{\<u,\xi\>\} Q_t(x,d\xi)
 =
\exp\{\<x,\psi(t,u)\>\}, \qquad x \in D, u \in U.
 \eeqlb
(The phrase ``homogeneous affine'' comes from the homogeneous
affine transformation $x \mapsto \<x,\psi(t,u)\>$.) We say the
HA-semigroup defined above is \textit{regular} if it is
stochastically continuous and the derivative
$(\partial\psi/\partial t) (0,u)$ exists for all $u\in U$ and is
continuous at $u=0$.

Clearly, the HA-semigroup satisfies the branching property
\eqref{2.1}, so the probability measure $Q_t(x,\cdot)$ is
infinitely divisible. To simplify the presentation, we assume that
$(Q_t)_{t\ge0}$ and all probabilities on $D$ possess finite first
absolute moments. Then the infinite divisibility of $Q_t(x,\cdot)$
and the special structure of $D$ imply that $\psi_2(t,u) =
\beta_{22}(t)u_2$ for some $\beta_{22}(t) \in \mbb{R}$ and
$\psi_1(t,u)$ has the representation
 \beqlb\label{4.2}
\psi_1(t,u) = \beta_{11}(t)u_1 + \beta_{12}(t)u_2 + \alpha(t)u_2^2
+ \int_D (e^{\<u,\xi\>} - 1 - u_2\xi_2) \mu(t,d\xi),
 \eeqlb
where $\alpha(t) \in \mbb{R}_+$, $(\beta_{11}(t), \beta_{12}(t))
\in D$ and $\mu(t,d\xi)$ is a $\sigma$-finite measure on $D$
supported by $D \setminus \{0\}$ such that
 \beqnn
\int_D \Big(|\xi_1| + |\xi_2|\land |\xi_2|^2\Big) \mu(t,d\xi) <
\infty;
 \eeqnn
see \cite{DL06}. From \eqref{4.2} and the semigroup property of
$(Q_t)_{t\ge0}$ it follows that
 \beqlb
\beta_{22}(r+t) &=& \beta_{22}(r)\beta_{22}(t), \label{4.3} \\
\beta_{11}(r+t) &=& \beta_{11}(r)\beta_{11}(t), \label{4.4} \\
\beta_{12}(r+t) &=& \beta_{11}(r)\beta_{12}(t) + \beta_{12}(r)
\beta_{22}(t), \label{4.5} \\
\alpha(r+t) &=& \beta_{11}(r)\alpha(t) + \alpha(r)\beta_{22}^2(t),
\label{4.6} \\
\mu(r+t,\cdot) &=& \int_D \mu(r,d\xi)Q_t(\xi,\cdot) +
\beta_{11}(r)\mu(t,\cdot) \label{4.7}
 \eeqlb
for any $r,t\ge0$.

The definition of SC-semigroups certainly applies to a
HA-semigroup. It was proved in Dawson and Li \cite{DL06} that if
$(\gamma_t)_{t\ge0}$ is a stochastically continuous SC-semigroup
associated with a regular HA-semigroup $(Q_t)_{t\ge0}$, then each
$\gamma_t$ is an infinitely divisible probability measure. Then we
have the representations
 \beqlb\label{4.8}
\int_{D} \exp\{\<u,\xi\>\} \gamma_t(d\xi)
 =
\exp\{\phi(t,u)\}, \qquad u\in U
 \eeqlb
and
 \beqlb\label{4.9}
\phi(t,u) = b_1(t)u_1 + b_2(t)u_2 + a(t)u_2^2 + \int_D
(e^{\<u,\xi\>} - 1 - u_2\xi_2) m(t,d\xi),
 \eeqlb
where $a(t) \in \mbb{R}_+$, $(b_1(t), b_2(t)) \in D$ and
$m(t,d\xi)$ is a $\sigma$-finite measure on $D$ supported by $D
\setminus \{0\}$ such that
 \beqnn
\int_D [|\xi_1| + |\xi_2|\land |\xi_2|^2] m(t,d\xi) < \infty.
 \eeqnn

\bproposition\label{p4.2} {\rm (\cite{DL06})} If
$(\gamma_t)_{t\ge0}$ is a stochastically continuous SC-semigroup
given by \eqref{4.8} and \eqref{4.9}, then for any $r,t\ge0$ we
have
 \beqlb
b_1(r+t) &=& b_1(r)\beta_{11}(t) + b_1(t), \label{4.10} \\
b_2(r+t) &=& b_1(r)\beta_{12}(t) + b_2(r) \beta_{22}(t) +
b_2(t)  \label{4.11} \\
a(r+t) &=& b_1(r)\alpha(t) + a(r)\beta_{22}^2(t) + a(t), \label{4.12} \\
m(r+t,\cdot) &=& \int_D m(r,d\xi)Q_t(\xi,\cdot) +
b_1(r)\mu(t,\cdot) + m(t,\cdot). \label{4.13}
 \eeqlb
\eproposition

The equations \eqref{4.10}--\eqref{4.13} give an alternative
expression of the property \eqref{2.4} and make it possible to
treat separately the coefficients in \eqref{4.9}. This leads to
some explicit analysis of the differentiability of $t\mapsto
\phi(t,u)$. In particular, if $\nu$ is an infinitely divisible
probability measure on $D$, we can define an SC-semigroup
$(\gamma_t)_{t\ge0}$ by \eqref{4.8} by letting
 \beqlb\label{4.14}
\phi(t,u) = \int_0^t \log\hat\nu(\psi(s,u)) ds, \qquad t\ge0, u\in
U,
 \eeqlb
where $\hat\nu$ is the characteristic function of $\nu$. In this
case, we call $(\gamma_t)_{t\ge0}$ a \textit{regular}
SC-semigroup. A simple but irregular SC-semigroup can be
constructed by letting $Q_t$ be the identity and letting $\gamma_t
= \delta_{(0,b_2(t))}$ where $b_2(t)$ is a discontinuous solution
of $b_2(r+t) = b_2(r) + b_2(t)$; see, e.g., \cite[p.37]{Sat99}.
This example shows that some condition on the function $t \mapsto
b_2(t)$ has to be imposed to get the regularity of the
SC-semigroup $(\gamma_t)_{t\ge0}$ given by \eqref{4.8} and
\eqref{4.9}. The proof of \cite[Theorem~3.1]{DL06} gives the
following

\btheorem\label{t4.1} {\rm (\cite{DL06})} Let $(\gamma_t)_{t\ge0}$
be a stochastically continuous SC-semigroup given by \eqref{4.8}
and \eqref{4.9}. Then the following conditions are equivalent:

\begin{itemize}

\item[(i)] $(\gamma_t)_{t\ge0}$ is regular;

\item[(ii)] $(\partial\phi/\partial t) (0,u)$ exists for every $u\in U$ and
is continuous at $u=0$;

\item[(iii)] $t\mapsto b_2(t)$ is absolutely continuous on $[0,\infty)$.

\end{itemize}\etheorem

Suppose that $(Q_t)_{t\ge0}$ is a HA-semigroup given by
\eqref{4.1} and $(\gamma_t)_{t\ge0}$ is an associated SC-semigroup
given by \eqref{4.8} and \eqref{4.9}. Let $P_t(x,\cdot) =
Q_t(x,\cdot) * \gamma_t(\cdot)$. Then $(P_t)_{t\ge0}$ is also a
Markov transition semigroup on $D$ and
 \beqlb\label{4.15}
\int_D \exp\{\<u,\xi\>\} P_t(x,d\xi) = \exp\{\<x,\psi(t,u)\> +
\phi(t,u)\}, \qquad x \in D, u \in U.
 \eeqlb

In general, a Markov transition semigroup on $D$ with
characteristic function of the form \eqref{4.15} is called an
\textit{affine semigroup}; see, e.g., \cite{DFS03}. We say the
affine semigroup is \textit{regular} if it is stochastically
continuous and the derivatives $(\partial\psi/\partial t) (0,u)$
and $(\partial\phi/\partial t) (0,u)$ exist for all $u\in U$ and
are continuous at $u=0$. Clearly, $(P_t)_{t\ge0}$ is regular if
and only if both $(Q_t)_{t\ge0}$ and $(\gamma_t)_{t\ge0}$ are
regular. Therefore, the above theorem gives a partial answer to
the open problem of characterizing all affine semigroups without
the regularity assumption; see \cite[Remark~2.11]{DFS03}. The
class of regular affine semigroups was characterized completely in
\cite{DFS03}. It was shown in \cite{DL06} that a regular affine
process arises naturally in a limit theorem for the difference of
a pair of reactant processes in a catalytic CBI-process.

\section{Measure-valued immigration processes}

\setcounter{equation}{0}

Let $E$ be a Lusin topological space, i.e., a homeomorph of a
Borel subset of a compact metric space. Recall that $M(E)$ is the
space of finite Borel measures on $E$ endowed with the topology of
weak convergence. A Markov process with state space $M(E)$ is
called an \textit{MB-process} if its transition semigroup
$(Q_t)_{t\ge0}$ satisfies the branching property \eqref{2.1}.
MB-processes appeared in Ji\v rina \cite{Jir58, Jir64} and
Watanabe \cite{Wat68} as high density limits of branching particle
systems. A very important special class of MB-processes, known as
Dawson-Watanabe processes, have been studied extensively in the
past decades; see, e.g., \cite{Daw93, Dyn02, Eth00, LeG99}. The
development of this subject has been stimulated from different
subjects including branching processes, interacting particle
systems, stochastic partial differential equations and non-linear
partial differential equations. The study of superprocesses has
also led to better understanding of results in those subjects.

Suppose that $(\gamma_t)_{t\ge0}$ is an SC-semigroup associated
with $(Q_t)_{t\ge0}$ and $(Q_t^\gamma)_{t\ge0}$ is defined by
\eqref{2.3}. A Markov process with transition semigroup
$(Q^\gamma_t)_{t\ge0}$ is called an \textit{immigration process}
associated with $(Q_t)_{t\ge0}$. This formulation of immigration
processes was given in Li \cite{Li95/6, Li96}. The intuitive
meaning of the immigration model is clear from the definition of
$(Q^\gamma_t)_{t\ge0}$. Clearly, this formulation essentially
includes all immigration mechanisms that are independent of the
inner population.

\btheorem\label{t5.1} {\rm (\cite{Li95/6})} A family of
probability measures $(\gamma_t)_{t\ge0}$ on $M(E)$ is an
SC-semigroup associated with $(Q_t)_{t\ge0}$ if and only if there
is an infinitely divisible probability entrance law $(K_t)_{t>0}$
for $(Q_t)_{t\ge0}$ such that
 \beqlb\label{5.1}
\log \int_{M(E)} e^{-\nu(f)} \gamma_t(d\nu)
 =
\int_0^t\bigg[\log \int_{M(E)} e^{-\nu(f)} K_s(d\nu)\bigg]ds,
\qquad t\ge 0, f\in B(E)^+.
 \eeqlb
\etheorem

Let us consider the case of a Dawson-Watanabe superprocess.
Suppose that $(P_t)_{t\ge0}$ is the transition semigroup of a
Borel right process $\xi$ with state space $E$ and
$\phi(\cdot,\cdot)$ is a branching mechanism given by
 \beqlb\label{5.2}
\phi(x,z) = b(x)z + c(x)z^2 + \int_0^\infty(e^{-zu}-1+zu) m(x,d
u), \qquad x\in E,\ z\ge0,
 \eeqlb
where $b\in B(E)$, $c\in B(E)^+$ and $(u\land u^2) m(x,du)$ is a
bounded kernel from $E$ to $(0,\infty)$. Then for each $f\in
B(E)^+$ the evolution equation
 \beqlb\label{5.3}
V_tf(x) + \int_0^tds\int_E \phi(y,V_sf(y)) P_{t-s}(x,dy) =
P_tf(x), \qquad t\ge0,\ x\in E
 \eeqlb
has a unique solution $V_tf\in B(E)^+$, and there is a Markov
semigroup $(Q_t)_{t\ge0}$ on $M(E)$ such that
 \beqlb\label{5.4}
\int_{M(E)} e^{-\nu(f)} Q_t(\mu,d\nu) = \exp\left\{-
\mu(V_tf)\right\}, \qquad f\in B(E)^+.
 \eeqlb
Clearly, $(Q_t)_{t\ge0}$ satisfies the branching property
\eqref{2.1}. A Markov process having transition semigroup
$(Q_t)_{t\ge0}$ is called a \textit{Dawson-Watanabe with
parameters $(\xi,\phi)$} or simply a
\textit{$(\xi,\phi)$-superprocess}. This process is a natural
generalization of the CB-process; see, e.g., \cite{Daw93}. The
family of operators $(V_t)_{t\ge0}$ form a semigroup, which is
called the \textit{cumulant semigroup} of the superprocess. Under
our hypotheses, $(Q_t)_{t\ge0}$ has a Borel right realization; see
\cite{Fit88, Fit92}.

Let $\mcr{K}(P)$ be the set of entrance laws $\kappa=
(\kappa_t)_{t>0}$ for the underlying semigroup $(P_t)_{t\ge0}$
that satisfy $\int_0^1\kappa_s(E)d s < \infty$. We endow
$\mcr{K}(P)$ with the $\sigma$-algebra generated by the mappings
$\kappa \mapsto \kappa_t(f)$ with $t>0$ and $f \in B(E)^+$. For
$\kappa\in \mcr{K}(P)$, set
 \beqlb\label{5.5}
S_t(\kappa,f) = \kappa_t(f) - \int_0^tds\int_E \phi(y,V_sf(y))
\kappa_{t-s}(dy), \qquad t>0,\ f\in B(E)^+.
 \eeqlb
Let $\mcr{K}^1(Q)$ denote the set of probability entrance laws $K
= (K_t)_{t>0}$ for the semigroup $(Q_t)_{t\ge0}$ of the
$(\xi,\phi)$-superprocess such that
 \beqnn
\int_0^1ds\int_{M(E)}\nu(E)K_s(d\nu) < \infty.
 \eeqnn

\btheorem\label{t5.2} {\rm (\cite{Li96})} A probability entrance
law $K\in\mcr{K}^1(Q)$ is infinitely divisible if and only if its
Laplace functional has the representation
 \beqlb\label{5.6}
\int_{M(E)} e^{-\nu(f)} K_t(d\nu) = \exp\bigg\{-S_t(\kappa,f)
-\int_{\mcr{K}(P)} \left(1-e^{-S_t(\eta,f)}\right) J(d\eta)
\bigg\},
 \eeqlb
where $\kappa \in \mcr{K}(P)$ and $J$ is a $\sigma$-finite measure
on $\mcr{K}(P)$ satisfying
 \beqnn
\int_0^1 d s\int_{\mcr{K}(P)} \eta_s(1) J(d\eta) <\infty.
 \eeqnn
\etheorem

The above theorem characterizes a class of infinitely divisible
probability entrance laws for $(Q_t)_{t\ge0}$. The right hand side
of \eqref{5.6} corresponds to an infinitely divisible probability
measure on $\mcr{K}(P)$. If $(\kappa_t)_{t>0}$ is given by
$\kappa_t = \mu P_t$, we have clearly $S_t(\kappa,f) = \mu(V_tf)$.
Let $\lambda \in M(E)$ and let $L$ be a $\sigma$-finite measure on
$M(E)$ satisfying
 \beqnn
\int_{M(E)} \nu(1) L(d\nu) <\infty.
 \eeqnn
We can define an infinitely divisible probability entrance law $K
\in \mcr{K}^1(Q)$ by
 \beqlb\label{5.7}
\int_{M(E)} e^{-\nu(f)} K_t(d\nu)
 =
\exp\bigg\{-\lambda(V_tf) - \int_{M(E)}
\left(1-e^{-\nu(V_tf)}\right) L(d\nu)\bigg\}.
 \eeqlb
This entrance law can be closed by an infinitely divisible
probability measure on $M(E)$. In this case, the transition
semigroup of the corresponding immigration process is given by
 \beqlb\label{5.8}
\int_{M(E)} e^{-\nu(f)} Q^\gamma_t(\mu,d\nu)
 &=&
\exp\bigg\{-\mu(V_tf) - \int_0^t\bigg[\lambda(V_sf) \nnm \\
 & &\qquad
+ \int_{M(E)} \Big(1-e^{-\nu(V_sf)}\Big) L(d\nu)\bigg]ds\bigg\}.
 \eeqlb
This is the case considered in Li \cite{Li92}. It was proved in Li
\cite{Li96} following the arguments of Fitzsimmons \cite{Fit88,
Fit92} that $(Q_t^\gamma)_{t\ge0}$ is a Borel right semigroup.

Needless to say, most of the theory of Dawson-Watanabe
superprocesses carries over to their associated immigration
processes and could be developed by techniques very close to those
of Dawson \cite{Daw93}. However, the immigration processes have
many additional structures, as might be expected from \eqref{2.3}
and \eqref{2.4}. A construction for the immigration processes were
given in Li \cite{Li02} by picking up measure-valued paths with
random times of birth and death. The construction was based on the
observation that any SC-semigroup is determined by a continuous
increasing measure-valued path $(\eta_t)_{t\ge0}$ and an entrance
rule $(G_t)_{t\ge 0}$. This structure yields a natural
decomposition of the immigration into two parts, the deterministic
part represented by $(\eta_t)_{t\ge0}$ and the random part
determined by $(G_t)_{t\ge0}$. The latter is an inhomogeneous
immigration process and can be constructed by summing up paths
$\{w_t: \alpha<t<\beta\}$ in the associated Kuznetsov process. By
analyzing the asymptotic behavior of the paths $\{w_t:
\alpha<t<\beta\}$ near the birth time $\alpha= \alpha(w)$, it was
shown in \cite{Li02} that almost all these paths start propagation
in an extension of the underlying space. Those combined with the
construction mentioned above give a full description of the
immigration phenomenon. As an application of the construction, Li
\cite{Li02} gave reformulations of some well-known results on
excessive measures in terms of stationary immigration
superprocesses. The immigration phenomena associated with
branching particle systems were studied in \cite{Li98}.

The state space of the immigration superprocess can be extended to
include some infinite measures; see, e.g., \cite{Li92, LW99}. With
such extensions, the immigration can be governed by a
$\sigma$-finite measure. A central limit theorem for the
$d$-dimensional super-Brownian motion with immigration was proved
in Li and Shiga \cite{LS95}, where the immigration is governed by
a deterministic $\sigma$-finite measure. When the governing
measure is the Lebesgue measure, the normalization function is
$t\sp{3/4}$ for $d=1$, $(t\log t)\sp{1/2}$ for $d=2$ and $t\sp
{1/2}$ for $d\geq 3$. The corresponding large deviation principle
was obtained in Zhang \cite{Zha04a} with the normalization
function $t$ in all dimensions and the speed function $t\sp {1/2}$
for $d=1$, $t/\log t$ for $d=2$ and $t$ for $d\geq 3$; see also
\cite{Zha05b}. The gap between the central limit theorem and the
large deviation principle was filled in Zhang \cite{Zha04b} by
establishing a moderate deviation principle. More precisely, she
proved that this immigration superprocess satisfies a large
deviations principle under the normalization $t^{1-\delta/4}$ for
$d=1$, $t^{1-\delta/2}(\log t)^{\delta/2}$ for $d=2$ and
$t^{1-\delta/2}$ for $d\geq 3$, where $\delta\in(0,1)$ is a
parameter; see also \cite{HL05, Zha05a}.

A super-Brownian motion with immigration governed by another
super-Brownian was introduced and studied in Hong and Li
\cite{HL99}. They established a central limit theorem for the
process which leads to Gaussian random fields in high dimensions.
For $d=3$ the field is spatially uniform, for $d\ge5$ its
covariance is given by the potential operator of the underlying
Brownian motion and for $d=4$ it involves a mixture of the two
kinds of fluctuations, which seems to be a new phenomenon in the
asymptotic behavior of measure-valued processes. There is a
similar phenomenon in the central limit theorem of the
corresponding occupation times obtained in Hong \cite{Hon02a} with
$d=6$ being critical. Some quenched mean limit theorems were
proved in Hong \cite{Hon05}. The moderate deviation principles for
the immigration superprocesses were established in \cite{Hon02b}.
Large deviation problems were studied in \cite{Hon03}, where the
speed functions are $t\sp{1/2}$ in $d=3$ and $t$ in $d\geq 4$. For
$d\neq 4$ the principle was accomplished by the well-known
G\"artner-Ellis theorem. In the critical dimension $d=4$, the
large deviation problem is much more difficult and only the limit
superior was established. We refer the reader to \cite{LW99} for a
more detailed survey on the early results measure-valued
immigration processes.

\section{Excursions and generalized immigration processes}

\setcounter{equation}{0}

Let $\alpha>0$ be a constant and $\{B(t): t\ge0\}$ a standard
Brownian motion. For any initial condition $x(0) = x \ge 0$ the
stochastic differential equation
 \beqlb\label{6.1}
dx(t) = \sqrt{2\alpha x(t)} d B(t), \qquad t\ge0
 \eeqlb
has a unique non-negative solution $\{x(t): t\ge0\}$, which is a
special case of the CB-process. This process is known as a
\textit{Feller branching diffusion} in the literature. The
transition semigroup $(Q_t)_{t\ge0}$ of the process is determined
by
 \beqlb\label{6.2}
\int_0^\infty e^{-z y}Q_t(x,dy) = \exp\{- xz(1+\alpha tz)^{-1}\},
\qquad t,x,z\ge0;
 \eeqlb
see, e.g., \cite[p.236]{IW89}. In view of the infinite
divisibility implied by \eqref{6.2}, there is a family of
\textit{canonical measures} $(\kappa_t)_{t>0}$ on $(0,\infty)$
such that
 \beqlb\label{6.3}
\int_0^\infty (1- e^{-zy})\kappa_t(dy) = z(1+\alpha tz/2)^{-1},
\qquad t>0,z\ge0.
 \eeqlb
Indeed, it is easy to check that
 \beqlb\label{6.4}
\kappa_t(dy) = (\alpha t)^{-2}e^{-y/\alpha t}dy, \qquad t,x>0.
 \eeqlb
Let $Q_t^\circ(x,dy)$ denote the restriction to $(0,\infty)$ of
the kernel $Q_t(x,dy)$. Since zero is a trap for the Feller
branching diffusion, $(Q_t^\circ) _{t\ge0}$ also constitute a
semigroup. Based on \eqref{6.2} and \eqref{6.3} it is not hard to
show that $\kappa_rQ_t^\circ = \kappa_{r+t}$ for all $r,t>0$. In
other words, $(\kappa_t) _{t>0}$ is an entrance law for
$(Q_t^\circ)_{t\ge0}$.

Let $W = C([0,\infty), \mbb{R}_+)$ and let $\tau_0(w) = \inf\{s>0:
w_s=0\}$ for $w\in W$. Let $W_0$ be the set of paths $w \in W$
such that $w_0 = w_t = 0$ for $t \ge \tau_0(w)$. We endow $W$ and
$W_0$ with the topology of locally uniform convergence. By the
theory of Markov processes, there is a unique $\sigma$-finite
measure $\mbf{Q}_\kappa$ on $(W_0,\mcr{B}(W_0))$ such that
 \beqlb\label{6.5}
\mbf{Q}_\kappa\{w_{t_1}\in dy_1,\cdots,w_{t_n}\in dy_n\} =
\kappa_{t_1}(dy_1)Q^\circ_{t_2-t_1}(y_1,dy_2)\cdots
Q^\circ_{t_n-t_{n-1}}(y_{n-1},dy_n)
 \eeqlb
for $0<t_1<t_2<\cdots<t_n$ and $y_1,y_2,\cdots,y_n\in (0,\infty)$;
see, e.g., \cite{PY82}. The measure $\mbf{Q}_\kappa$ is known as
the \textit{excursion law} of the Feller branching diffusion.
Roughly speaking, \eqref{6.5} asserts that $\{w_t: t>0\}$ under
$\mbf{Q}_\kappa$ is a Feller branching diffusion with
one-dimensional distributions $\{\kappa_t: t>0\}$. The Feller
branching diffusion can be reconstructed from the excursion law
$\mbf{Q}_\kappa$ in the following way: Fix $x\ge0$ and let $N(dw)$
be a Poisson random measure on $W_0$ with intensity $x
\mbf{Q}_\kappa (dw)$. Let $x(0)=x$ and
 \beqlb\label{6.6}
x(t) = \int_{W_0} w_t N(dw), \qquad t>0.
 \eeqlb
Then $\{x(t): t\ge0\}$ is a weak solution of \eqref{6.1}; see
\cite[Theorem~4.1]{PY82}.

Let $b(\cdot)$ be a non-negative and locally Lipschitz function on
$\mbb{R}_+$ satisfying the linear growth condition. A non-negative
diffusion process $\{y(t): t\ge0\}$ can be defined by the
stochastic differential equation
 \beqlb\label{6.7}
dy(t) = \sqrt{2\alpha y(t)} d B(t) + b(y(t))dt, \qquad t\ge0.
 \eeqlb
This process can be constructed from a Feller branching diffusion
and a Poisson random measure based on $\mbf{Q}_\kappa(dw)$ as
follows. Let $\{x(t): t\ge0\}$ be a Feller branching diffusion and
let $N(ds,du,dw)$ be a Poisson random measure on $\mbb{R}_+^2
\times W_0$ with intensity $dsdu\mbf{Q}_\kappa(dw)$. We assume
that $\{x(t): t\ge0\}$ and $\{N(ds,du,dw)\}$ are independent.

\bproposition\label{p6.1} {\rm (\cite{FL04})} There is a unique
strong solution of the stochastic equation
 \beqlb\label{6.8}
y(t) = x(t) + \int_0^t\int_0^{b(y(s))}\int_{W_0} w_{t-s}
N(ds,du,dw),  \qquad t\ge 0.
 \eeqlb
Moreover, the solution $\{y(t): t\ge0\}$ of the above equation is
a weak solution of \eqref{6.7}. \eproposition

This proposition is a consequence of Fu and Li
\cite[Theorem~4.1]{FL04}, where more general results on
measure-valued processes were given. In particular, if $b(x)
\equiv b$ is a constant, $\{y(t): t\ge0\}$ is a CBI-process
associated with the Feller branching diffusion; see \cite{PY82}.
In the general case, we may regard $\{y(t): t\ge0\}$ as a
\textit{generalized CBI-process}.

The approach of stochastic equations driven by Poisson random
measures based on the excursion law has more substantial
applications in constructions of some measure-valued diffusions.
Let us look at an example of this type involving a stochastic
flow. Suppose that $h$ is a continuously differentiable function
on $\mbb{R}$ such that both $h$ and $h^\prime$ are
square-integrable. Then the function
 \beqlb\label{6.9}
\rho(x) = \int_{\mbb{R}}h(y-x)h(y) dy, \qquad x\in\mbb{R}
 \eeqlb
is twice continuously differentiable with bounded derivatives
$\rho^\prime$ and $\rho^{\prime\prime}$. Let $m$ be a
$\sigma$-finite Borel measure on $\mbb{R}$ and $q(\cdot,\cdot)$ a
non-negative Borel function on $M(\mbb{R}) \times \mbb{R}$ such
that there is a constant $K$ such that
 \beqlb\label{6.10}
\int_{\mbb{R}} q(\mu,y) m(dy) \le K(1 + \|\mu\|), \qquad \mu\in
M(\mbb{R}),
 \eeqlb
and for each $R>0$ there is a constant $K_R>0$ such that
 \beqlb\label{6.11}
\int_{\mbb{R}} |q(\mu,y) - q(\nu,y)| m(dy) \le K_R \|\mu-\nu\|
 \eeqlb
for $\mu$ and $\nu\in M(\mbb{R})$ satisfying $\mu(\mbb{R}) \le R$
and $\nu(\mbb{R}) \le R$, where $\|\cdot\|$ denotes the total
variation of the signed measure. Let us consider the following
martingale problem of an $M(\mbb{R})$-valued process $\{Y_t:
t\ge0\}$: For each $\phi\in C^2(\mbb{R})$,
 \beqlb\label{6.12}
M_t(\phi) := Y_t(\phi) - Y_0(\phi) - \rho(0)\int_0^t
Y_s(\phi^{\prime\prime}) ds - \int_0^tds\int_{\mbb{R}}
\phi(y)q(Y_s,y) m(dy)
 \eeqlb
is a continuous martingale with quadratic variation process
 \beqlb\label{6.13}
\<M(\phi)\>_t = 2\alpha\int_0^t Y_s(\phi^2) ds + \int_0^t
ds\int_{\mbb{R}^2} \rho(x-y) \phi^\prime(x)\phi^\prime(y)
Y_s^2(dx,dy).
 \eeqlb

Let $W(dt,dy)$ be a time-space white noise on $[0,\infty) \times
\mbb{R}$ based on the Lebesgue measure; see, e.g., \cite{Wal86}.
By \cite[Lemma~3.1]{DLW01} or \cite[Lemma~1.3]{Wan97}, for any
$r\ge0$ and $a\in \mbb{R}$ the stochastic equation
 \beqlb\label{6.14}
x(t) = a + \int_r^t\int_{\mbb{R}} h(y-x(s))W(ds,dy), \qquad t\ge r
 \eeqlb
has a unique continuous solution $\{x(r,a,t): t\ge r\}$, which is
a Brownian motion with quadratic variation $\rho(0)dt$. Indeed,
the system $\{x(r,a,t): t\ge r; a\in\mbb{R}\}$ determines an
isotropic stochastic flow. Fix $\mu \in M(\mbb{R})$ and let
$N_0(da,dw)$ be a Poisson random measure on $\mbb{R} \times W_0$
with intensity $\mu(da)\mbf{Q}_\kappa(dw)$ and $N(ds,da,du,dw)$ a
Poisson random measure on $\mbb{R}_+ \times \mbb{R} \times
\mbb{R}_+ \times W_0$ with intensity $ds m(da)du
\mbf{Q}_\kappa(dw)$. Suppose that $\{W(dt,dy)\}$, $\{N_0(da,dw)\}$
and $\{N(ds,da,du,dw)\}$ are independent of each other.

\btheorem\label{t6.1} {\rm (\cite{DL03})} There is a unique strong
solution of the stochastic equation
 \beqlb\label{6.15}
Y_t
 &=&
\int_{\mbb{R}}\int_{W_0} w(t) \delta_{x(0,a,t)} N_0(da,dw) \nnm \\
 & &
+ \int_0^t\int_{\mbb{R}}\int_0^{q(Y_s,a)}\int_{W_0} w(t-s)
\delta_{x(s,a,t)} N(ds,da,du,dw), \qquad t>0.
 \eeqlb
Furthermore, if we set $Y_0 = \mu$, the process $\{Y_t: t\ge0\}$
is a measure-valued diffusion process solving the martingale
problem given by \eqref{6.12} and \eqref{6.13}. \etheorem

In view of \eqref{6.15}, we may regard $\{Y_t: t\ge0\}$ as a
generalized immigration superprocess carried by the stochastic
flow given by \eqref{6.14}. The stochastic equation \eqref{6.15}
is substantial for the construction of this measure-valued
diffusion process, for the uniqueness of solution of the
martingale problem given by \eqref{6.12} and \eqref{6.13} still
remains open.

\section{Generalized Mehler semigroups}

\setcounter{equation}{0}

Let $H$ be a real separable Hilbert space and let $(T_t)_{t\ge0}$
be a strongly continuous semigroup of linear operators on $H$. A
family of probability measures $(\gamma_t)_{t\ge0}$ on $H$ is
called an \textit{SC-semigroup} associated with $(T_t)_{t\ge0}$ if
it satisfies
 \beqlb\label{7.1}
\gamma_{r+t} = (\gamma_r \circ T_t^{-1}) * \gamma_t, \qquad
r,t\ge0.
 \eeqlb
This is clearly the special case of \eqref{2.4} with
$Q_t(x,\cdot)\equiv \delta_{T_tx}$. If \eqref{7.1} is satisfied,
we can define a Markov transition semigroup $(Q^\gamma_t)_{t\ge0}$
on $H$ by
 \beqlb\label{7.2}
Q^\gamma_tf(x) = \int_H f(T_tx+y)\mu_t(dy), \qquad x\in H,f\in
B(H),
 \eeqlb
which is called a \textit{generalized Mehler semigroup} associated
with $(T_t)_{t\ge0}$. The corresponding Markov process is a
generalized OU-process; see \cite{BRS96}.

According to a result of Schmuland and Sun \cite{SS01}, if
$(\gamma_t)_{t\ge0}$ is a solution of \eqref{7.1}, each $\gamma_t$
is an infinitely divisible probability measure. By Linde
\cite[p.75 and p.84]{Lin86}, we have the following representation
of the characteristic functional:
 \beqlb\label{7.3}
\hat\gamma_t(a)
 &=&
\exp\bigg\{i\<b_t,a\> - \frac{1}{2}\<R_ta,a\>  \nnm \\
 & &\qquad
+ \int_{H^\circ} \Big(e^{i\< x,a\>} -1 - i\< x,a\>
1_{[0,1]}(\|x\|)\Big) M_t(dx)\bigg\}, \qquad t\ge0, a\in H,
 \eeqlb
where $b_t\in H$, $R_t$ is a symmetric, positive-definite nuclear
operator on $H$, and $M_t$ is a $\sigma$-finite measure (L\'evy
measure) on $H^\circ := H \setminus \{0\}$ satisfying
 \beqnn
\int_{H^\circ} (1\land\|x\|^2) M_t(dx) < \infty.
 \eeqnn

\btheorem\label{t7.1} {\rm (\cite{DLSS04})} Suppose that
$(\gamma_t)_{t\ge0}$ is a family of probability measures on $H$.
If there is a family of infinitely divisible probabilities
$(\nu_s)_{s>0}$ such that $\nu_{r+t} = \nu_r \circ T_t^{-1}$ for
all $r,t>0$ and
 \beqlb\label{7.4}
\hat\gamma_t(a) = \exp\bigg\{\int_0^t \log\hat\nu_s(a) ds\bigg\},
\qquad t\ge0, a\in H,
 \eeqlb
then $(\gamma_t)_{t\ge0}$ is an SC-semigroup. Conversely, if
$(\gamma_t)_{t\ge0}$ is an SC-semigroup given by \eqref{7.3} and
if $t\mapsto b_t$ is absolutely continuous, then the
characteristic functional $\hat\gamma_t$ has representation
\eqref{7.4}. \etheorem

The above theorem gives a characterization for the generalized
Mehler semigroup $(Q^\gamma_t)_{t\ge0}$. If there is an infinitely
divisible probability measure $\nu_0$ on $H$ such that $\nu_t =
\nu_0\circ T_t^{-1}$, we say that $(\gamma_t)_{t\ge0}$ and
$(Q^\gamma_t)_{t\ge0}$ are \textit{regular}. In this case, the
function $t \mapsto \hat \gamma_t(a)$ is differentiable for every
$a\in H$. It was proved in Bogachev \textit{et al.}\ \cite{BRS96}
that a cylindrical Gaussian SC-semigroup satisfying this
differentiability condition can be extended into a real Gaussian
SC-semigroup in an enlargement of $H$ and the corresponding
OU-process can be constructed as the strong solution to a
stochastic differential equation. Those results were extended to
the general non-Gaussian case in \cite{FR00}. A simple and nice
necessary and sufficient condition for the differentiability of $t
\mapsto \hat \gamma_t(a)$ was given by van Neerven \cite{Nee00}.
Some powerful inequalities for regular generalized Mehler
semigroups were proved in R\"ockner and Wang \cite{RW03} and Wang
\cite{Wan05a}.

It was observed in Dawson \textit{et al.}\ \cite{DLSS04} that the
OU-processes corresponding to an irregular generalized Mehler
semigroup usually have no right continuous realizations. Under the
second moment assumption, Dawson and Li \cite{DL04} studied the
construction of OU-processes corresponding to centered but
irregular SC-semigroups. Based on Theorem~\ref{t7.1} they showed
that each centered SC-semigroup is uniquely determined by an
infinitely divisible probability measure on the entrance space
$\tilde H$ for the semigroup $(T_t)_{t\ge0}$, which is an
enlargement of $H$. They proved that a centered SC-semigroup can
always be extended to a regular one on the entrance space. Those
results provide an approach to the study of irregular generalized
Mehler semigroups with which one can reduce some of their analysis
to the framework of \cite{BRS96, FR00, Wan05a}.

\section{Fluctuation limits of immigration processes}

\setcounter{equation}{0}

Fluctuation limits of branching particle systems and
superprocesses have been studied extensively. Since those systems
are usually unstable, in many cases one uses time-dependent
scalings which lead to time-inhomogeneous OU-processes; see, e.g.,
\cite{BG86, DFG89} and the references therein. For subcritical
branching systems with immigration, it is usually easy to find a
stationary distribution. In the study of fluctuation limits of
those systems, we can use a time-independent scaling, which lead
to homogeneous OU-processes. Fluctuation limits of this kind were
studied in \cite{GL98, GL00, Li99, Li00b, LZ05}.

Let $A_0$ be the generator of a conservative Feller transition
semigroup on $\mbb{R}^d$ such that $C^2(\mbb{R}^d) \subseteq
\mcr{D}(A_0)$ and $A_0f \in C(\mbb{R}^d)$ for every $f\in
C^2(\mbb{R}^d)$. We fix a strictly positive function $b(\cdot) \in
C(\mbb{R}^d)^+$ which is bounded away from zero. Let
$(P_t)_{t\ge0}$ be the semigroup generated by $A := A_0-b$ and let
$\phi_0$ be a continuous function given by
 \beqlb\label{8.1}
\phi_0(x,z) = c(x)z^2+\int_0^\infty (e^{-zu}-1+zu)n(x,du), \qquad
x\in \mbb{R}^d, z\geq 0,
 \eeqlb
where $c(\cdot) \in C(\mbb{R}^d)^+$ and $u^2n(x,du)$ is a bounded
kernel from $\mbb{R}^d$ to $(0,\infty)$. Then the evolution
equation
 \beqlb\label{8.2}
V_tf(x) + \int_0^tds\int_{\mbb{R}^d} \phi_0(y,V_sf(y))
P_{t-s}(x,dy) = P_tf(x), \qquad t\geq 0, x\in \mbb{R}^d
 \eeqlb
defines a cumulant semigroup $(V_t)_{t\geq 0}$. Given $m \in
M(\mbb{R}^d)$, we can define the transition semigroup
$(Q_t^m)_{t\ge0}$ of an immigration superprocess with state space
$M(\mbb{R}^d)$ by
 \beqlb\label{8.3}
\int_{M(\mbb{R}^d)} e^{-\nu(f)}Q_t^m(\mu,d\nu)
 =
\exp\bigg\{-\mu(V_tf)-\int_0^t m(V_sf)ds\bigg\}, \qquad f \in
C(\mbb{R}^d)^+.
 \eeqlb
Since $(P_t)_{t\geq 0}$ is a Feller semigroup, the immigration
superprocess has a Hunt realization. In particular, it has a
c\`adl\`ag realization; see, e.g., \cite[p.221]{Sha88}.

It is easy to see that $Q_t^m(\mu, \cdot)$ converges as $t\to
\infty$ to the probability measure $Q_\infty(\cdot)$ on
$M(\mbb{R}^d)$ given by
 \beqlb\label{8.4}
\int_{M(\mbb{R}^d)} e^{-\nu(f)}Q_\infty^m(d\nu)
 =
\exp\bigg\{-\int_0^\infty m(V_sf)ds\bigg\},
 \qquad f \in C(\mbb{R}^d)^+.
 \eeqlb
Clearly, $Q_\infty(\cdot)$ is the unique equilibrium of the
semigroup $(Q_t^m)_{t\ge0}$. Moreover, we have
 \beqlb\label{8.5}
\int_{M(\mbb{R}^d)} \nu(f) Q_\infty^m(d\nu) = \lambda(f), \qquad f
\in C(\mbb{R}^d)^+,
 \eeqlb
where $\lambda \in M(\mbb{R}^d)$ is defined by
 $$
\lambda = \int_0^\infty m P_s ds.
 $$

It is a natural problem to investigate the asymptotic fluctuation
of the immigration superprocess around the long-term average
$\lambda$ as the branching mechanism $\phi_0$ decreases to zero. A
result of this type is formulated as follows. For any integer
$k\ge1$ let $\phi_k(x,z) = \phi_0(x,z/k)$. Then $\phi_k(x,z) \to
0$ as $k\to \infty$. Suppose that $\{Y_t^{(k)}: t\geq 0\}$ is a
c\`adl\`ag immigration superprocess with parameters $(A,\phi_k,m)$
and $Y_0^{(k)} = \lambda$. Let
 \beqlb\label{8.6}
Z_t^{(k)} = k[Y_t^{(k)}-\lambda], \qquad t\geq 0.
 \eeqlb
Let $\mcr{S}(\mbb{R}^d)$ denote the Schwartz space of rapidly
decreasing functions on $\mbb{R}^d$. That is, each $f\in
\mcr{S}(\mbb{R}^d)$ is belong to $C^\infty(\mbb{R}^d)$ and for
each integer $n\ge 1$ and each non-negative integer-valued vector
$\alpha = (\alpha_1,\cdots,\alpha_d)$ we have
 $$
\lim_{|x|\to \infty}|x|^n|\partial^\alpha f(x)|=0,
 $$
where
 $$
\partial^\alpha f(x) = \frac{\partial^{|\alpha|}}
{\partial x_1^{\alpha_1} \cdots x_d^{\alpha_d}} f(x_1,\cdots,x_d)
 $$
and $|\alpha| = \alpha_1 + \cdots + \alpha_d$. The topology of
$\mcr{S}(\mbb{R}^d)$ is defined by the sequence of semi-norms
 $$
f\mapsto p_n(f) := \sup\{(1+|x|^n)|\partial^\alpha f(x)|: x\in
\mbb{R}^d, |\alpha| \le n\}, \qquad n=0,1,2,\cdots.
 $$
Let $\mcr{S}^{\prime}(\mbb{R}^d)$ denote the dual space of
$\mcr{S}(\mbb{R}^d)$ equipped with the strong topology. Then both
$\mcr{S}(\mbb{R}^d)$ and $\mcr{S}^{\prime}(\mbb{R}^d)$ are nuclear
spaces; see, e.g., \cite[p.107]{Sch80}. It is easy to see that
$\{Z_t^{(k)}: t\geq 0\}$ has sample paths in $D([0,\infty),
\mcr{S}^\prime(\mbb{R}^d))$.

\btheorem\label{t8.1} {\rm (\cite{GL98})} As $k\to \infty$, the
finite dimensional distributions of $\{Z_t^{(k)}: t\geq 0\}$
converge to those of the $\mcr{S}^{\prime} (\mbb{R}^d)$-valued
Markov process $\{Z_t:t\geq 0\}$ with $Z_0=0$ and with transition
semigroup $(T_t)_{t\geq 0}$ defined by
 \beqlb\label{8.7}
\int_{\mcr{S}^{\prime}(\mbb{R}^d)} e^{i\<\nu,f\>} T_t(\mu,d\nu)
 =
\exp\bigg\{i\<\mu,P_tf\> + \int_0^t \lambda(\phi_0(-iP_sf))
ds\bigg\}, \qquad f \in \mcr{S}(\mbb{R}^d),
 \eeqlb
where $\phi_0(-iP_sf)$ is given by \eqref{8.1} with $z$ replaced
by $-iP_sf(x)$.\etheorem

The above theorem was improved in \cite{LZ05}, where it was proved
that $\{Z_t^{(k)}: t\ge 0\}$ converges to $\{Z_t: t\ge 0\}$ weakly
in $D([0,\infty), \mcr{S}^\prime(\mbb{R}^d))$. Indeed, the
fluctuation limit theorem was also formulated in \cite{LZ05} in a
suitable Sobolev space. For any integer $n\ge0$ we define the
Sobolev space
 $$
H^n(\mbb{R}^d) = \{f\in \mcr{S}^{\prime}(\mbb{R}^d):
\partial^\alpha f \in L^2(\mbb{R}^d) \mbox{~whenever~} |\alpha|\le n\}
 $$
with the norm $\|\cdot\|_n$ defined by
 \beqnn
\|f\|_n^2 = \sum_{|\alpha|\le n}\int_{\mbb{R}^d}
|\partial^\alpha{f}(x)|^2 dx.
 \eeqnn
Let $H^{-n}(\mbb{R}^d)$ be the strong topological dual of
$H^n(\mbb{R}^d)$. It is well-known that $H^{-n}(\mbb{R}^d)$ can be
identified as a subspace of $\mcr{S}^\prime(\mbb{R}^d)$ and
 \beqlb\label{8.8}
\mcr{S}^{\prime}(\mbb{R}^d) \supseteq H^m(\mbb{R}^d) \supseteq
H^n(\mbb{R}^d) \supseteq \mcr{S}(\mbb{R}^d)
 \eeqlb
for any integers $m\le n$ with continuous embeddings; see, e.g.,
\cite[Theorem~5.5]{BN73}. Now we have

\btheorem\label{t8.2} {\rm (\cite{LZ05})} For any integer $n>d+2$
the process $\{Z_t: t\ge 0\}$ has a realization in $D([0,\infty),
H^{-n}(\mbb{R}^d))$ and $\{Z_t^{(k)}: t\ge 0\}$ converges weakly
to $\{Z_t: t\ge 0\}$ in $D([0,\infty), H^{-n}(\mbb{R}^d))$.
\etheorem

By the above theorem, $\{Z_t: t\ge 0\}$ is a generalized
OU-process in the real separable Hilbert space
$H^{-n}(\mbb{R}^d)$. This puts the process into the framework of
generalized Mehler semigroup of the last section and makes it
possible to derive regularities and properties of the processes
from the existing literature; see, e.g., \cite{BRS96, DLSS04,
FR00, RW03, Wan05a}.

The limiting generalized OU-process obtained in above can live in
a much smaller state space. Let us consider the case where $A_0 =
\Delta$ and $\phi_0(x,z) = c(x)z^2/2$. In this case, the
corresponding generalized OU-process solves the Langevin equation
 \beqlb\label{8.9}
dZ_t = dW_t + \Delta Z_tdt - bZ_tdt, \qquad t\ge0,
 \eeqlb
where $\{W_t: t\ge 0\}$ is a time-space white noise with intensity
$c(x)dt\lambda(dx)$; see, e.g., \cite{Li99}. Given $Z_0$ the
solution of \eqref{8.9} is represented by
 \beqlb\label{8.10}
Z_t = Z_0P_t + \int_0^t\int_{\mbb{R}^d}p_{t-s}(x,\cdot)W(ds,dx),
\qquad t\ge0,
 \eeqlb
where $p_t(x,\cdot)$ denotes the density of $P_t(x,\cdot)$. If
$d=1$, the process $\{Z_t: t\ge 0\}$ has a version in
$L^2(\mbb{R})$. Indeed, it is well-known that $Z_0P_t \in
L^2(\mbb{R})$ whenever $Z_0 \in L^2(\mbb{R})$. On the other hand,
we have
 \beqnn
& &\mbf{E}\bigg[\int_{\mbb{R}}\bigg(\int_0^t\int_{\mbb{R}}
p_{t-s}(x,y) W(ds,dx)\bigg)^2dy\bigg]  \\
 &=&
\int_{\mbb{R}}dy\int_0^tds\int_{\mbb{R}} p_{t-s}(x,y)^2
c(x)\lambda(dx)  \\
 &\le&
\int_0^t\frac{1}{\sqrt{2\pi(t-s)}}ds\int_{\mbb{R}}
c(x)\lambda(dx)  \\
 &<&
\infty.
 \eeqnn
Then the second term on the right hand side of \eqref{8.10} exists
almost surely in $L^2(\mbb{R})$. It follows that
 \beqlb\label{8.11}
\int_{L^2(\mbb{R})} e^{i\<\nu,f\>} T_t(\mu,d\nu)
 =
\exp\bigg\{i\<\mu,P_tf\> - \int_0^t \lambda(c|P_sf|^2) ds\bigg\},
\qquad f\in L^2(\mbb{R})
 \eeqlb
defines a generalized Mehler semigroup $(T_t)_{t\geq 0}$ on
$L^2(\mbb{R})$. This semigroup is clearly irregular on the state
space $L^2(\mbb{R})$, but the characteristic functional of the
corresponding SC-semigroup is differentiable in time.
Measure-valued catalysts for superprocess were introduced by
Dawson and Fleischmann \cite{DF91}. One may also study fluctuation
limits of immigration superprocesses with measure-valued
catalysts. In such case the resulting SC-semigroup may have
non-differentiable characteristic functionals; see \cite{DLSS04}.

\medskip

{\small

\textbf{Acknowledgements.} This work was supported by the Creative
Research Group Fund of the National Natural Science Foundation of
China under grant number 10121101. I own special thanks to
Professors M.F.\,Chen and F.Y.\,Wang who suggested that I write
this survey. I would also like to thank Professors W.M.\,Hong and
M.\,Zhang for helpful discussions on related topics.}

\end{document}